\documentclass[10pt,reqno]{amsart}
\usepackage{graphicx}
\usepackage{amssymb,amsfonts,epsfig,xspace,color,enumerate}
\usepackage{verbatim}
\usepackage{amsmath}

\newtheorem{thm}{Theorem}[section]
\newtheorem{cor}[thm]{Corollary}

\newtheorem{prop}[thm]{Proposition}
\theoremstyle{definition}
\newtheorem{defn}[thm]{Definition}
\newtheorem{rem}[thm]{Remark}
\newtheorem{exa}[thm]{Example}

\newcommand{\sinc}{\textnormal{sinc}\xspace}
\newcommand{\supp}{\textnormal{supp}\,}
\newcommand{\closedspan}{\overline{\textnormal{span}}}
\newcommand{\Real}{\mathbb{R}}
\newcommand{\sgn}{\textnormal{sgn}\,}
\newcommand{\Z}{\mathbb{Z}}

\newcommand{\R}{\mathbb{R}}

\newcommand{\spa}{\textnormal{span}\,}

\newcommand{\lam}{\lambda}

\begin{document}


\title{Irregular and multi--channel sampling of operators}

\author{Yoon Mi Hong}
\address{ Department of Mathematical Sciences, KAIST,
Daejeon, 305-701, Republic of Korea} 

\author{G\"otz E. Pfander} \address{
School of Engineering and Science,
Jacobs University, 28759 Bremen, Germany}
\email{ymhong@amath.kaist.ac.kr, g.pfander@jacobs-university.de}

%


\begin{abstract}
The classical sampling theorem for bandlimited functions has
recently been generalized to apply to so-called bandlimited
operators, that is, to operators with band-limited Kohn-Nirenberg
symbols. Here, we discuss operator sampling versions of two of the
most central extensions to the classical sampling theorem. In
irregular operator sampling, the sampling set is not periodic with
uniform distance. In multi-channel operator sampling, we obtain complete
information on an operator by multiple operator sampling outputs.
\end{abstract}


\maketitle

\section{Introduction}
The so-called classical sampling theory addresses the problem of
determining and reconstructing functions
from countably many values that are attained on a discrete subset of
the real line. The fundamental result in this theory is attributed to
Whittaker, Kotel'nikov and Shannon. It asserts that a function
bandlimited to an interval of length $\Omega$ can be recovered from
the values of the function sampled regularly  at $\Omega$ values per
unit interval.

During the last few years the herein considered sampling theory for
operators has been developed. It is motivated by the operator
identification problem in communications engineering. There, the
objective is to identify a channel operator from knowledge of the
channel's action on a  chosen input signal. A well known identification result states, for
example, that any time--invariant channel operator is fully
determined by its action on the Dirac impulse. Already in the 1960s,
Kailath \cite{Kai} and Bello \cite{Bel} proclaimed that this simple
identifiability result on time--invariant operators could be
generalized to slowly time--varying operators, that is, to operators
whose spreading functions are supported on sets of measure less than
or equal to one. The spreading function is the symplectic Fourier
transform of the operator's Kohn-Nirenberg
symbol. The assertions of Kailath and Bello were confirmed in
\cite{KoPf, PfWa06}.

The identifiability results in \cite{KoPf, PfWa06} use weighted sums
of regularly spaced delta impulses as identifiers for
Hilbert--Schmidt operators with bandlimited Kohn--Nirenberg symbols.
The discrete support of such a tempered distribution we shall refer
to in the following as sampling set for operator sampling. Together
with the fact that the classical sampling theorem can be seen as a
special case of the identifiability results for bandlimited
operators \cite{KoPf, PfWa06,Pf} --- consider a bandlimited function as a
multiplication operator which can be determined from its action on a regularly spaced sum of
Dirac impulses --- this has led to the development of the herein
considered sampling theory for operators.

In this paper, we state and prove operator sampling versions of key
generalizations of the classical sampling theorem for functions.

One central extension of the classical sampling theorem considers irregular or nonuniform sampling sets, see \cite{Mar} and references therein. In practice, sampling with uniform distance is hardly realizable because of imperfections in data acquisition devices or perturbations during the collecting of data. As it is similarly challenging to generate regularly spaced sums of  impulses for operator sampling, the consideration of irregularly placed impulses as identifiers for operators is natural. Here, we will give separation and Beurling density
results for operator sampling which resemble corresponding results in the classical sampling theory for functions.

The feasibility of the classical sampling theory for very large bandwidth signals is limited by the sampling rates achievable in state of the art hardware. This problem is addressed through multi-channel sampling as pioneered by Papoulis \cite{Pap}. Multi--channel sampling employs a number of samplers in parallel,
thereby allowing the acquisition of samples to be carried out in
each channel at a fraction of the sampling rate foreseen by the
classical sampling theorem. Multi--channel sampling in the theory of
operator sampling employs similarly the combination of multiple
outputs from sampling procedures in order to reduce the rate at
which impulses have to be sent. In addition, multi--channel
sampling for operators, that is, multiple output sampling for
operators, allows to identify operators whose Kohn--Nirenberg
symbols are only bandlimited to an area of measure larger than one, as shown below.
Note that in single output operator sampling, only operators with
bandlimitations given by sets of measure less than one can be
identified. Larger bandwidth of the Kohn--Nirenberg symbol cannot be
compensated by an increase of the so-called sampling rate.
\footnote{Operator sampling is not simply an higher dimensional
analogue of the 1-d Shannon sampling theorem. In the  case of an
operator acting on $L^2(\R)$, the operator's  2--dimensional
Kohn--Nirenberg symbol is to be determined from a signal defined on
$\R$. No access to sample values of the Kohn--Nirenberg symbol is
given, as is the case in 2-dimensional Shannon sampling theory.}

We formulate and prove our results for Hilbert--Schmidt operators.
Means for generalizing such results to non compact operators are outlined in
\cite{PfWa, Pf}.

This paper is organized as follows. In Section 2, we give some
background for operator sampling and introduce operator Paley-Wiener
spaces. In Section 3 we state the uniform operator sampling
result for operator Paley-Wiener spaces as given in \cite{PfWa07}. We
include a new proof of this result, a proof that allows for
generalizations to the setting of irregular and multi-channel
operator sampling. We give a
generalization to operator classes which have not necessarily
bandlimited Kohn--Nirenberg symbols. In Section 4, we develop
irregular operator sampling for operator Paley-Wiener spaces. Also, we
consider irregular sampling of operators whose Kohn-Nirenberg symbol is not bandlimited
in view of Kramer's Lemma setting. Multi-channel operator sampling
is discussed in Section 5.

\section{Preliminaries}
The Fourier transform on $L^2(\Real^d)$ is densely defined by
$$\mathcal{F}(f)(\xi)=\hat{f}(\xi)=\int_{\mathbb{R}^d}f(t)\,e^{-2\pi i t\cdot\xi}dt,\quad f\in L^1(\Real^d)\cap L^2(\Real^d).$$
 Similarly, the symplectic Fourier transform
$\mathcal{F}_s: L^2(\mathbb{R}^{2d})\longrightarrow
L^2(\mathbb{R}^{2d})$ is given by
$$\mathcal{F}_sf(t,\nu)=\int_{\mathbb{R}^{d}}\int_{\mathbb{R}^{d}}f(x,\xi)\,e^{-2\pi i(\nu\cdot x-\xi\cdot t)}dxd\xi,\quad a.e. \ t,\nu\in\R.$$
For brevity of notation, we shall  refrain from marking equalities and inequalities that hold in the $L^2$-sense with the customary $a.e.$ whenever the context is unambiguous.

Let
$S_\Omega=\prod_{k=1}^d[-\frac{\Omega_k}{2},\frac{\Omega_k}{2}]$
with $\Omega=(\Omega_1,\Omega_2,\cdots,\Omega_d)$.
\begin{defn}
The Paley-Wiener space with bandwidth ${\Omega}$ is defined by
$$PW(S_\Omega)=\left\{f\in L^2(\mathbb{R}^d)\,:\,\supp \hat{f}\subseteq
S_\Omega\right\}.$$
\end{defn}
It is known that if an expansion for $f\in PW(S_\Omega)$ converges
in the norm of $PW(S_\Omega)$, then it converges pointwise and
uniformly over $\Real^d$. Note that $PW(S_\Omega)$ is isometrically
isomorphic to $L^2(S_\Omega)$ due to Plancherel's
theorem. 
For simplicity of notation, we shall denote by $L^2(S_\Omega)$ the
subspace of $L^2(\Real^{d})$ which consists of $L^2(\Real^{d})$ functions
supported on $S_\Omega$.


 The classical sampling theorem shown independently by
Whittaker, Kotel'nikov and Shannon is generalized by the following
{\it oversampling} theorem. It is based on possibly collecting
samples  more often than the sampling rate
prescribes. The sampling rate, usually called
Nyquist-Landau rate, for a function in $PW(S_\Omega)$ is defined to
be the Lebesgue measure of the set $S_\Omega$.

Here and in the following we use the notation $$A(F)\asymp
B(F),\quad F\in\mathcal F,$$ if there exist positive constants $c$
and $C$ such that $c A(F)\leq B(F)\leq C A(F)$ for all objects $F$ in  the set $
\mathcal F$.

\begin{thm}\label{WKSthm}
For $\Omega, T>0$ with $T_k\Omega_k<1$ for all $k$ and
$\varphi\in PW(S_{\frac{2}{T}-\Omega})$ 
with $\hat{\varphi}=1$ on $S_{\Omega}$, we have a sampling expansion
\begin{equation}\label{wks}
f(t)=T\sum_{n\in\mathbb{Z}^d}f(nT)\varphi(t-nT),\quad f\in
PW(S_\Omega).
\end{equation}
Moreover, (\ref{wks}) is stable in the sense that
$$\|f\|^2\asymp \sum_{n\in\mathbb{Z}^d} |f(nT)|^2,\quad f\in PW(S_\Omega).$$
\end{thm}
In practice, a stable sampling expansion guarantees that
perturbations in the sampling output and in the reconstruction procedure
are controlled by error bounds on the input function and vice versa. 

The development of operator sampling necessitates the use of some
rudimentary distribution theory. The space of distributions chosen here
is the dual of the Feichtinger algebra  $S_0(\Real^d)$. The dual  $S_0'(\Real^d)$ is a Banach space with $\mathcal S (\Real^d)\subsetneq S_0(\Real^d)\subsetneq L^2(\Real^d)\subsetneq S_0'(\Real^d)\subsetneq \mathcal{S}'(\Real^d)$, where  $\mathcal{S}$ denotes the Schwartz class of rapidly decaying functions and
$\mathcal{S}'$ its dual, the space of tempered distribution. 
There are several equivalent definitions of the Feichtinger algebra $S_0$
\cite{FeZi}. We choose the characterization of $S_0$ via the
short time Fourier transform.

\begin{defn}
The Feichtinger algebra is defined by
$$S_0(\mathbb{R}^d)=\{f\in L^2(\mathbb{R}^d)\,:\,V_g f(t,\nu)\in L^1(\mathbb{R}^{2d})\},$$
where $V_g f(t,\nu)=\langle f,M_\nu T_t g\rangle$ is the short-time
Fourier transform of $f$ with respect to the Gaussian
$g(x)=e^{-\pi\|x\|^2}$. The norm on $S_0(\R^d)$ is given by
$\|f\|_{S_0}=\|V_g f\|_{L^1}$.
\end{defn}

In this paper we consider the sampling
problem for Hilbert-Schmidt operators only.
\begin{defn}
The class of  Hilbert-Schmidt operators $HS(L^2(\R^d))$ consists of bounded linear operators on
$L^2(\mathbb{R}^d)$ which can be represented as  integral operators of the form
\begin{equation*}
Hf(x)=\int\kappa_H(x,t)f(t)dt, \quad f\in L^2(\R^d),
\end{equation*}
with kernel $\kappa_H\in L^2(\mathbb{R}^{2d})$.
\end{defn}

The linear space of Hilbert-Schmidt operators
$HS(L^2(\mathbb{R}^{d}))$ becomes a Hilbert space if it is endowed
with the Hilbert space structure of $L^2(\mathbb{R}^{d})$, that is,  by
$$
    \langle H_1,H_2\rangle_{HS}=\langle
    \kappa_{H_1},\kappa_{H_2}\rangle_{L^2}.
$$
In view of pseudodifferential operators, the Kohn-Nirenberg symbol
$\sigma_H$ \cite{Fol, KoNi} of a Hilbert-Schmidt operator $H$ is
given by
$$\sigma_H(x,\xi)=\int \kappa_H(x,x-t)\,e^{-2\pi it\cdot\xi}dt.$$
It leads to the operator representation
$$Hf(x)=\int\sigma_H(x,\xi)\hat{f}(\xi)\,e^{2\pi ix\cdot\xi}d\xi, \quad f\in L^2(\R^d).$$
In time-frequency analysis and communication engineering, the
spreading function $\eta_H$ of a Hilbert-Schmidt operator $H$ is
commonly considered. It is given by
$$\eta_H(t,\nu)=\int\kappa_H(x,x-t)\,e^{-2\pi ix\cdot\nu}dx,$$
and leads to
\begin{equation}\label{weakHS}
Hf(x)=\int\int\eta_H(t,\nu)M_\nu T_t f(x)\,dtd\nu, \quad f\in L^2(\R^d),
\end{equation}
where the time-shift (translation) and frequency-shift (modulation)
operators $T_t$ and $M_\nu$ are defined by $T_t f(x)=f(x-t)$ and
$\widehat{M_\nu f}(\gamma)=\hat{f}(\gamma-\nu)$, respectively. That
is, a Hilbert-Schmidt operator $H$ is a continuous superposition of
translation and modulation operators with coefficient function
$\eta_H$. The identity (\ref{weakHS}) is understood weakly, namely
$$\langle Hf,g\rangle=\int\int \eta_H(t,\nu) \langle M_\nu T_t
f,g\rangle\, dtd\nu, \quad g\in L^2(\Real^d).$$ As
$\eta_H=\mathcal{F}_s\sigma_H$, operators with band-limited
Kohn-Nirenberg symbols are operators whose spreading functions are
compactly supported.

In communications, the time-varying operator $H$ is also commonly
represented by its time-varying impulse response $h_H(t,x)$ with
$$Hf(x)=\int h_H(t,x)f(x-t)dt,$$
where $h_H(t,x)=\kappa_H(x,x-t)=\int \eta_H(t,\nu)\,e^{2\pi
ix\nu}d\nu$ a.e.
 Note that
$$\|H\|_{HS}=\|\kappa_H\|_{L^2}=\|h_H\|_{L^2}=\|\sigma_H\|_{L^2}=\|\eta_H\|_{L^2}.$$


\begin{defn}
The operator Paley-Wiener
space of operators bandlimited to $S\subseteq \R^{2d}$ is
$$OPW(S)=\{H\in HS(L^2(\mathbb{R}^d))\,:\,\supp \mathcal F_s \sigma_H\subseteq S\}.$$
\end{defn}
In the literature, operators with $\supp \mathcal F_s \sigma_H\subseteq [a_1,b_1]{\times}\ldots{\times} [a_{2d},b_{2d}]$ are commonly referred to as underspread or slowly time--varying operators if \\ volume$([a_1,b_1]{\times}\ldots{\times} [a_{2d},b_{2d}]=(b_1-a_1)\cdot \ldots \cdot (b_{2d}-a_{2d})\leq 1$ and as overspread operators else (see \cite{KoPf} and references within).

We formulate the operator identification and sampling problems as follows.
\footnote{See \cite{KoPf} for a more general concept of operator
identification.}

\begin{defn}
An operator class $\mathcal H \subseteq HS(L^2(\R^d))$ is
identifiable if all $H\in \mathcal H$ extend to a domain containing
a so--called identifier $f\in S_0'(\R^d)$ with
 \begin{equation}\label{iden}
\|H\|_{HS}\asymp \|Hf\|_{L^2},\quad H\in\mathcal{H}.
\end{equation}

The operator class $\mathcal H \subseteq HS(L^2(\R^d))$ permits operator sampling if one can choose $f$ in  (\ref{iden}) with discrete support in $\R^d$ in the distributional sense. In that case, $\supp f$ is called sampling set for $\mathcal H$.
\end{defn}

Note that $H\in OPW(S)$ with $S$ compact can be extended to a
bounded linear operator $H:S_0'(\Real^d)\longrightarrow
L^2(\Real^d)$
\cite{CoGr,PfWa}. This implies that proving identifiability of a
Hilbert-Schmidt operator by an element in $S_0'$ is equivalent to
providing the lower bound $A$ in (\ref{iden}), as an upper bound is
given by $B=\|f\|_{S_0'}$ since
 $$\|Hf\|_{L^2}
\leq \|H\|_\textrm{op}\|f\|_{S_0'} \leq \|H\|_{HS} \|f\|_{S_0'},\quad H\in OPW(S).$$


Given a separable Hilbert space $X$, a sequence of elements
$\{f_k\}_{k\in\mathbb{Z}}$ in $X$ is called a frame for $X$ if
$$\sum_{k\in\mathbb{Z}}|\langle f,f_k\rangle|^2\asymp\|f\|_X^2,\quad  f\in X.$$
To each frame $\{f_k\}_{k\in\mathbb{Z}}$ for $X$ exists a so-called dual frame
$\{\widetilde f_k\}_{k\in\mathbb{Z}}$  of $\{f_k\}_{k\in\mathbb{Z}}$ for $X$ with
\begin{equation*}
  f=\sum_{k\in\Z} \langle f, f_k\rangle \widetilde f_k =\sum_{k\in\Z} \langle f, \widetilde f_k\rangle  f_k, \quad f\in X. 
\end{equation*}
Moreover, a frame which does not form a
frame if we remove any element from it is called a Riesz basis, or, also, exact frame. A sequence
$\{f_k\}_{k\in\Z}$ is called a Riesz sequence if it is a Riesz basis
for $\closedspan\{f_k\}_{k\in\Z}$ \cite{Chr, Gro01, KoCh}.

\section{Uniform sampling of Hilbert-Schmidt
operators}\label{section:uniform}


Theorem \ref{WKSthm} states that a bandlimited, square integrable
function can be reconstructed by its values sampled at a
sufficiently dense sampling grid. In this section, we first consider
the sampling problem for operators whose Kohn-Nirenberg symbols are
bandlimited in the time-frequency plane. We shall state and prove
all theorems for $d=1$ for convenience.

Operators with rectangular bandlimitation on their Kohn-Nirenberg
symbols are the starting point of operator sampling \cite{PfWa07,Pf},

\begin{thm}\label{uniform}
For $\Omega,T,T'>0$ and $0<T'\Omega \leq{ T\Omega \leq 1}$, choose
$\varphi\in
PW([-(\frac{1}{T}-\frac{\Omega}{2}),\frac{1}{T}-\frac{\Omega}{2}])$
with $\hat{\varphi}=1$ on $[-\frac{\Omega}{2},\frac{\Omega}{2}]$ and
$r\in L^\infty(\Real)$ with $\supp r\subset [-T+T',T]$ and $r=1$ on
$[0,T']$. Then  $OPW([0,T']{\times}[-\displaystyle
\tfrac{\Omega}{2},\tfrac{\Omega}{2}])$ permits operator sampling as
\begin{equation*}
\|H\|_{HS}=\sqrt{T}\, \|H\sum_{n\in\Z}\delta_{nT}\|_{L^2},\quad H\in OPW([0,T']{\times}[-\displaystyle
\tfrac{\Omega}{2},\tfrac{\Omega}{2}]),
\end{equation*}
 and operator reconstruction is possible by means of
 $$h_H(t,x)=r(t)T\sum_{n\in\mathbb{Z}}(H\sum_{k\in\Z}\delta_{kT})(t+nT)\varphi(x-t-nT).$$
\end{thm}

In \cite{KoPf}, the proof of the identifiability of
$OPW([0,T']{\times}[-\tfrac{\Omega}{2},\tfrac{\Omega}{2}])$ is based
on the unitarity of the Zak transform. For clarity and to indicate
directions for generalizations of this theorem, we prove
Theorem~\ref{uniform} through elementary orthonormal basis
expansions based on Fourier series.

\begin{proof}
%
For almost every $t\in\R$, we have $\eta_H(t,\cdot)\in
L^2[-\tfrac{\Omega}{2},\tfrac{\Omega}{2}]\subseteq
L^2[-\tfrac{1}{2T},\tfrac{1}{2T}]$ and, by expanding $\eta_H(t,\nu)$
with respect to the orthonormal basis\\
$
\{\sqrt T\, e^{-2\pi i(t+nT)\nu}\}_{n\in\Z}$ of
$L^2[-\tfrac{1}{2T},\tfrac{1}{2T}]$, we obtain
\begin{eqnarray}
\eta_H(t,\nu)&=&\sum_{n\in\Z}\langle\eta_H(t,\nu),\sqrt T\,e^{-2\pi
i(t+nT)\nu}\rangle \sqrt T\,e^{-2\pi i(t+nT)\nu}\notag\\
&=& T\sum_{n\in\Z} h_H(t,t+nT)\,e^{-2\pi i(t+nT)\nu}, \quad
\nu\in[-\tfrac{\Omega}{2},\tfrac{\Omega}{2}],\label{ONBexp}
\end{eqnarray}
so
$$\eta_H(t,\nu)=T\hat{\varphi}(\nu)\sum_{n\in\Z} h_H(t,t+nT)\,e^{-2\pi i(t+nT)\nu}, \quad \nu\in\Real,$$
with $\varphi$ chosen to satisfy $\varphi\in
PW([-(\frac{1}{T}-\frac{\Omega}{2}),\frac{1}{T}-\frac{\Omega}{2}])$
and $\hat{\varphi}=1$ on $[-\frac{\Omega}{2},\frac{\Omega}{2}]$.
Then, for almost every $t\in \R$
\begin{eqnarray}
h_H(t,x)&=&T\sum_{n\in\Z} h_H(t,t+nT)\int\hat{\varphi}(\nu)e^{-2\pi
i(t+nT)\nu}e^{2\pi i\nu x}d\nu\notag\\
&=&T\sum_{n\in\Z} h_H(t,t+nT)\,\varphi(x-t-nT),\quad x\in\R.\label{h_H}
\end{eqnarray}
 On the other hand, we have $(H\sum_k\delta_{kT})(x)=\sum_k
h_H(x-kT,x)\in L^2(\R)$ so that $(H\sum_{k}\delta_{kT})(t+nT)=\sum_k
h_H(t+nT-kT,t+nT)$. Now, $r(t)(H\sum_{k\in\Z}\delta_{kT})(t+nT)=h_H(t,t+nT)$
for $t\in[0,T']$ since $\supp h_H(\cdot,x)\subseteq
[0,T']\subseteq[0,T]$. With (\ref{h_H}), this gives
\begin{equation}\label{expw/id}
h_H(t,x)=r(t)T \sum_{n\in\Z}
(H\sum_{k\in\Z}\delta_{kT})(t+nT)\varphi(x-t-nT),\quad x\in\R, \ \text{a.e. } t\in\R,
\end{equation}
a formula which contains the identifier $\sum_{k\in\Z}\delta_{kT}$ for
$OPW([0,T']{\times}[-\tfrac{\Omega}{2},\tfrac{\Omega}{2}])$. Note that
the series in (\ref{expw/id}) converges pointwise and uniformly over
$\Real$ in $x$.

Moreover,
since (\ref{ONBexp}) was an orthonormal basis expansion,
Parseval's identity gives
$$\|\eta_H(t,\cdot)\|_{L^2}^2{=}T\sum_{n\in\Z}|h_H(t,t+nT)|^2{=}T\sum_{n\in\Z}|r(t)(H\sum_{k\in\Z}\delta_{kT})(t+nT)|^2,
\ \ \text{a.e. } t\in\R$$
so that
$$\|h_H(t,\cdot)\|_{L^2}^2=T\sum_{n\in\Z}|r(t)(H\sum_{k\in\Z}\delta_{kT})(t+nT)|^2,\quad \text{a.e. } t\in\R,$$
and
\begin{eqnarray*}
\|h_H\|_{L^2}^2&=&\int_0^{T'}\|h_H(t,\cdot)\|_{L^2}^2dt
\ =\ {T}\int_0^{T'}\sum_{n\in\Z}|h_H(t,t+nT)|^2dt\\
&=&{T}\int_0^{T}\sum_{n\in\Z}|h_H(t,t+nT)|^2dt
\ =\ T\int_{0}^{T}\sum_{n\in\Z}|(H\sum_{k\in\Z}\delta_{kT})(t+nT)|^2dt\\
&=&T\int_\Real|(H\sum_{k\in\Z}\delta_{kT})(t)|^2dt
\ =\ T\|H\sum_{k\in\Z}\delta_{kT}\|_{L^2}^2
\end{eqnarray*}
as $h_H(t,t+nT)$ and $(H\sum_k\delta_{kT})(t+nT)$ vanish on the
interval $[T',T]$. The operator class
$OPW([0,T']{\times}[-\displaystyle
\tfrac{\Omega}{2},\tfrac{\Omega}{2}])$ is identifiable by
$\sum_{k\in\Z}\delta_{kT}$ as we showed
$\|H\|_{HS}=\|h_H\|_{L^2}=\sqrt{T}\|H\sum_{k\in\Z}\delta_{kT}\|_{L^2}$.
\end{proof}

One crucial ingredient in the proof above is the fact that for each
$t$, the set $\mathcal E_t=\{\sqrt T\, e^{-2\pi i(t+nT)\nu}\}_{n\in\Z}$ is an orthonormal
basis for $L^2[-\frac 1 {2T},\frac 1 {2T}]$. Note that the
 functionals corresponding to $\mathcal E_t$ depend on $t$ which is necessary to
associate $h_H(t,t+nT)$ with $Hg(t+nT)$ for some identifier $g$.
Another important ingredient is the support condition on
$h_H(\cdot,x)$. It guarantees that no aliasing in the infinite
summation $H(\sum_{k\in\Z}\delta_{kT})(t+nT)$ takes place as for each $t$ the
sum is reduced to a single non zero summand $h_H(t,t+nT)$.

We assume in Theorem \ref{uniform}
that the area of the rectangle $[0,
T']{\times}[-\frac{\Omega}{2},\frac{\Omega}{2}]$ is less than or
equal to 1. This assumption coincides with the one given in Kailath's conjecture for
identifiability of such operator classes \cite{Kai}. For
$1<T'\Omega\leq T\Omega$,  perfect reconstruction of $h_H(t,x)$ from
its samples is not possible since the sampling rate $\frac{1}{T}$ is
strictly less than the Nyquist-Landau rate $\Omega$ for $h_H(t,\cdot)$. In
this case, not only operator sampling, but also operator
identification by any tempered distribution as single input signal
is not possible as shown in \cite{Pf08}.

Now we extend Theorem \ref{uniform} to the case where $h_H(t,\cdot)$ lies in a  shift-invariant space other than the Paley-Wiener space. 
Given a Riesz sequence
$\{\varphi(\cdot-nT)\}_{n\in\Z}$ in $L^2(\Real)$, let
$$V_T(\varphi)=\overline{\spa}\{\varphi(\cdot-nT)\}=\big\{\sum_{n\in\Z}c_n\varphi(\cdot-nT)\,:\,\{c_n\}_{n\in\Z}\in l^2\big\}.$$
Let $\mathcal{H}_{T,\varphi}\subseteq HS(L^2(\Real))$ consist of
integral operators $H$ 
with $h_H\in L^2[0,T]\otimes
V_T(\varphi)$.
Note that different from the operator Paley-Wiener setup, not   each such operator maps boundedly   $S_0'$ to $L^2$.

We require the shift-invariant space $V_T(\varphi)$ to be a
reproducing kernel Hilbert space \cite{Hig}. 

\begin{defn}
A Hilbert space $X$ of complex-valued functions on a given domain
$D\neq\emptyset$ is a reproducing kernel Hilbert space if there
exists a kernel $k(s,t)$ defined on $D\times D$ satisfying
 $k(\cdot,t)\in X$ for all $t\in D$ and
 $f(t)=\langle f(\cdot),k(\cdot,t)\rangle_X$ for all $f\in X$ and $ t\in
D$. Such a function $k(s,t)$ is called a reproducing kernel.
\end{defn}

For example, if $\varphi$ is a complex-valued integrable function
well-defined everywhere in $\Real$ and satisfies
$$\sum_{n\in\Z}|\varphi(t+n)|^2<\infty, \quad t\in [0,1],$$
then $V_1(\varphi)$ is a reproducing kernel Hilbert space \cite{KK}.
 Alternatively, if $\varphi$ is
continuous and belongs to the Wiener amalgam space $W(L^\infty, l^1)$,
that is, to the subspace of $L^2(\R)$ defined by the norm
$$\|\varphi\|_{W(L^\infty, l^1)}=\sum_{n\in\Z}\sup_{t\in [0,1]}|\varphi(t+n)|<\infty,$$
then $V_1(\varphi)$ is a reproducing kernel Hilbert space as well
\cite{AlGr}.
\begin{thm}
Assume that $V_T(\varphi)$ is a reproducing kernel Hilbert space and
its reproducing kernel $k(s,t)$ satisfies the condition that
$\{k(\cdot,t+nT)\}_{n\in\Z}$ is a frame for $V_T(\varphi)$ for each
$t\in [0,T ]$. Then $\sum_{k\in\Z}\delta_{kT}$ identifies
$\mathcal{H}_{T,\varphi}$, that is,
$$\|H\|_{HS}\asymp\|H\sum_{k\in\Z}\delta_{kT}\|_{L^2},\quad H\in\mathcal{H}_{T,\varphi}.$$
The reconstruction of operators is possible by
$$h_H(t,x)=\chi_{[0,T]}(t)\sum_{n\in\Z}(H\sum_{k\in\Z}\delta_{kT})(t+nT)k_n^*(x,t),$$
where $\{k_n^*(\cdot,t)\}_{n\in\Z}$ is a dual frame of
$\{k(\cdot,t+nT)\}_{n\in\Z}$ for a.e. $t\in [0,T]$.
\end{thm}
\begin{proof}
Since $h_H(t,\cdot)\in V_T(\varphi)$ and
$\{k(\cdot,t+nT)\}_{n\in\Z}$ is a frame for $V_T(\varphi)$ for each
$t\in [0,T]$,
$$h_H(t,x)=\sum_{n\in\Z}\langle h_H(t,\cdot),k(\cdot,t+nT)\rangle k_n^*(x,t)$$
where $\{k_n^*(\cdot,t)\}_{n\in\Z}$ is a dual frame of
$\{k(\cdot,t+nT)\}_{n\in\Z}$. Since $k(s,t)$ is a reproducing kernel
of $V_T(\varphi)$,
$$h_H(t,x)=\sum_{n\in\Z} h_H(t,t+nT)k_n^*(x,t)$$
and
\begin{equation}\label{EqGP1}
  \|h_H(t,\cdot)\|_{L^2}^2\asymp \sum_{n\in\Z} |h_H(t,t+nT)|^2.
\end{equation}
Formally, we have
$(H\sum_{k\in\Z}\delta_{kT})(t+nT)=\sum_{k\in\Z} h_H(t+nT-kT,t+nT)$
so that
$$r(t)(H\sum_{k\in\Z}\delta_{kT})(t+nT)=h_H(t,t+nT)$$
where $r(t)=\chi_{[0,T]}(t)$. 
Together with (\ref{EqGP1}), we note $H\sum_{k\in\Z}\delta_{kT}\in L^2(\Real)$ and we conclude
\begin{eqnarray*}
\|h_H\|_{L^2}^2=\|H\|_{HS}^2&\asymp&\int_0^T\sum_{n\in\Z}|(H\sum_{k\in\Z}\delta_{kT})(t+nT)|^2dt
=\|H\sum_{k\in\Z}\delta_{kT}\|_{L^2}^2.
\end{eqnarray*}
\end{proof}
\begin{exa}\rm
Let $T=1$ and take $\varphi(t)=\chi_{[0,1)}(t)$. Then $V_1(\varphi)$
is a reproducing kernel Hilbert space and its reproducing kernel
$k(s,t)$ allows $\{k(\cdot,t+n)\}_{n\in\Z}$ to be a frame, in fact,
an orthonormal basis for $V_1(\varphi)$ since
$k(s,t)=\sum_{n\in\Z}\chi_{[n,n+1)}(s)\chi_{[n,n+1)}(t)$. Hence, for
$t\in [0,1]$,
$$k(s,t+n)=\sum_{m\in\Z} \chi_{[m,m+1)}(s)\chi_{[m,m+1)}(t+n)=\sum_{m\in\Z}
\chi_{[m,m+1)}(s)\delta_{n,m}=\chi_{[n,n+1)}(s).$$ We conclude
that $\sum_k\delta_k$ identifies $\mathcal{H}=\{H\in
HS(L^2(\Real))\,:\,h_H\in L^2[0,1]\otimes V_1(\varphi)\}$. For
$H\in\mathcal H$, the
 kernel $\kappa_H(x,t)$ is a step function along diagonals
$x=t+c,\, c\in\Real$.
\end{exa}

\section{Irregular sampling of Hilbert-Schmidt
operators}\label{section:irregular}
 First, we provide the background on
irregular sampling of functions that is needed to develop operator
sampling results based on irregular sampling sets.

\subsection{Irregular Sampling in Paley-Wiener spaces}
\begin{defn}\label{def_density}
Let $\Lambda=\{\lam_k\}_{k\in\Z}\subseteq\Real$, with
$\lam_k<\lam_{k+1}$, $k\in\Z$.
\begin{enumerate}
    \item $\Lambda$ is a set of sampling, also referred to as stable sampling
    set, for \\ $PW([-\frac \Omega 2,\frac \Omega 2])$ if
    $$\|f\|_{L^2}^2\asymp \sum_{k\in\Z}|f(\lam_k)|^2,\quad f\in PW([-\tfrac
    \Omega 2,\tfrac \Omega 2]).$$
    \item $\Lambda$ is a set of interpolation for $PW([-\frac \Omega 2,\frac \Omega
    2])$ if the interpolation or moment problem
    $$f(\lam_k)=c_k,\quad k\in\Z,$$
    has a solution in $PW([-\frac \Omega 2,\frac \Omega 2])$ for every
    $\{c_k\}\in l^2(\Z)$.
    \item $\Lambda$ is uniformly discrete if,
    $$(\lam_{k+1}-\lam_k)\geq\delta>0, \quad k\in\Z.$$
    In this case $\delta$ is called a separation constant.
    \item $\Lambda$ is relatively uniformly discrete if $\Lambda$ is
    a finite union of uniformly discrete sets.
    \item The upper and lower Beurling densities are  defined, respectively, by
    $$D^+(\Lambda)=\limsup_{h\rightarrow\infty}\frac{n^+(h)}{h}\quad\textrm{and}\quad D^-(\Lambda)=\liminf_{h\rightarrow\infty}\frac{n^-(h)}{h},$$
    where for $h>0$, $n^+(h)$ and $n^-(h)$ are the largest number and smallest number of points from $\Lambda$ in
    $[x-\frac h 2,x+\frac h 2)$, $x\in\Real$, respectively. If
    $D^+(\Lambda)=D^-(\Lambda)$, then we say that $\Lambda$ has uniform
    Beurling density $D(\Lambda)=D^+(\Lambda)=D^-(\Lambda)$.

\end{enumerate}
\end{defn}
We recall necessary and sufficient conditions on the Beurling density of
a set $\Lambda=\{\lam_k\}_{k\in\Z}$ for its nonharmonic sequence to be a
frame or a Riesz sequence for $L^2[-\tfrac{\Omega}{2},\tfrac{\Omega}{2}]$ \cite{Sei}.
\begin{thm}~\label{Sei}
\begin{enumerate}
    \item $\Lambda$ is a set of sampling for
    $PW([-\tfrac{\Omega}{2},\tfrac{\Omega}{2}])$ if and only if
    $\{e^{-2\pi i\lam_k\xi}\}_{k\in\Z}$ is a frame for
    $L^2[-\tfrac{\Omega}{2},\tfrac{\Omega}{2}]$. Moreover, for $\Lambda$ being relatively uniformly discrete, a necessary
    condition for
    $\{e^{-2\pi i\lam_k\xi}\}_{k\in\Z}$ to be a frame for
    $L^2[-\tfrac{\Omega}{2},\tfrac{\Omega}{2}]$ is $D^-(\Lambda)\geq
    \Omega$, and a sufficient condition is $D^-(\Lambda)>\Omega$.
    \item $\Lambda$ is a set of interpolation for
    $PW([-\tfrac{\Omega}{2},\tfrac{\Omega}{2}])$ if and only if
    $\{e^{-2\pi i\lam_k\xi}\}_{k\in\Z}$ is a Riesz sequence 
    in
    $L^2[-\tfrac{\Omega}{2},\tfrac{\Omega}{2}]$. Moreover, for $\Lambda$ being uniformly discrete, a necessary
    condition for
    $\{e^{-2\pi i\lam_k\xi}\}_{k\in\Z}$ to be a Riesz sequence in
    $L^2[-\tfrac{\Omega}{2},\tfrac{\Omega}{2}]$ is $D^+(\Lambda)\leq
    \Omega$, and a sufficient condition is $D^+(\Lambda)<\Omega$.
\end{enumerate}
\end{thm}

In general, it is highly non-trivial to determine whether a set
$\{e^{-2\pi i\lam_k\xi}\}_{k\in\Z}$ is a Riesz basis for $L^2[-\frac
\Omega 2,\frac \Omega 2]$. A famous affirmative result was given by
Kadec \cite{Kadec}.

\begin{thm}(Kadec's 1/4-theorem)\label{Kad}
Let $\Lambda=\{\lam_k\}_{k\in\Z}\subset\Real$. If there is
$L\geq 0$ such that
\begin{equation}
  |\lam_k-\tfrac k \Omega|\leq L<\tfrac 1 {4\Omega}, \quad k\in\Z,  \label{equation:Kadec}
\end{equation}
then $\{e^{-2\pi i\lam_k\xi}\}_{k\in\Z}$ is a Riesz basis for
$L^2[-\frac \Omega 2,\frac \Omega 2]$, and $\tfrac 1 {4\Omega}$ is the best
possible constant for (\ref{equation:Kadec}) to hold.
\end{thm}

\subsection{Irregular sampling in Operator Paley-Wiener spaces}

\begin{defn}
A sequence $\Lambda=\{\lam_k\}_{k\in\Z}$ in $\Real$ is a set of
sampling for an operator class $\mathcal{H}$, if for some
$\{c_k\}_{k\in\Z}\in l^\infty(\Z)$, we have $\sum_{k\in\Z} c_k
\delta_{\lam_k}\in S_0'(\R)$ and $\sum_{k\in\Z} c_k \delta_{\lam_k}$
identifies $\mathcal{H}$.
\end{defn}

Operator sampling is operator identification with discretely
supported identifiers. Consequently, irregular operator sampling for
$OPW([0,T]{\times}[-\frac \Omega 2,\frac \Omega 2])$ is {\it
a-priori}  possible only if $T\Omega \leq 1$ \cite{KoPf,PfWa06}.

\begin{thm}\label{thmGP}
  $T\Omega\leq 1$ is a necessary condition for the existence of a sampling set for $OPW([0,T]{\times}[-\frac \Omega 2,\frac \Omega 2])$.
\end{thm}

Analogous to Theorem~\ref{Sei}, we have the following result.

\begin{thm}\label{thmformerprop}
If $\Lambda=\{\lam_k\}_{k\in\Z}$, $\lam_{k+1}> \lam_k$, is uniformly
discrete, then a necessary condition for $\Lambda$ being a set of
sampling for $OPW([0,T]{\times}[-\frac \Omega 2,\frac \Omega 2])$ is
$D^-(\Lambda)\geq \Omega$ and a sufficient condition is
$D^-(\Lambda)>\Omega$ and $\lam_{k+1}-\lam_k\geq
T$.
\end{thm}
\begin{proof}
Assume that $\Lambda$ is a set of sampling for
$OPW([0,T]{\times}[-\frac \Omega 2,\frac \Omega 2])$ with
$D^{-}(\Lambda)<\Omega$. Then there exist $C$ and $\tilde{C}>0$ such
that
$$\|H\|_{HS}^2 {\leq} C \| \sum_{k\in\Z} \kappa_H(\,\cdot\, ,\lambda_k) \|^2 {\leq} \tilde{C} \sum_{k\in\Z}  \| \kappa_H(\,\cdot\,,\lambda_k) \|^2,\ \  H\in OPW([0,T]{\times}[-\tfrac \Omega 2,\tfrac \Omega 2])  , $$
since for each $x$ there exist at most $\lfloor\frac T \delta\rfloor
+1$ nonzero summands above where $\delta=\inf_k(\lam_{k+1}-\lam_k)$.
But as $D^- (\Lambda) <\Omega$,  there exists $f\in PW([-\frac
\Omega 2,\frac \Omega 2])$ with $\| f \| =1$ and $\sum_k|
f(\lambda_k)|^2\leq 1 /(2\tilde{C})$. Defining $H$ by
$\kappa_H(x,y)=f(y)$ for $0\leq x-y \leq T$, $\kappa_H=0$ else, we
have $\kappa_H(x,\lambda_k)= f(\lambda_k) $ for $ \lambda_k \leq x
\leq \lambda_k +T$ and $0$ else.

We have
$$ T= \|H\|_{HS}^2 \leq \tilde C \sum_{k\in\Z}\|\kappa_H(x,\lambda_k)\|^2 =\tilde C \sum_{k\in\Z}T|f(\lambda_k)|^2
\leq\tilde C T / 2\tilde C = T / 2, $$ a
contradiction.

Now we shall establish the  sufficient condition for $\Lambda$ to be a set of sampling. As $\Lambda$ is uniformly discrete and
$D^-(\Lambda)>\Omega$, $\{e^{-2\pi i \lam_k\xi}\}_{k\in\Z}$ is a
frame for $L^2[-\frac \Omega 2,\frac \Omega 2]$.  Theorem~\ref{Sei}
(1) implies that $\Lambda$ is a set of sampling for $PW([-\frac \Omega
2,\frac \Omega 2])$ and for $t\in\R$,
$$\|h_H(t,\cdot)\|_{L^2}^2\asymp\sum_{n\in\Z}|h_H(t,t+\lam_n)|^2,  \quad H\in OPW([0,T]{\times}[-\tfrac \Omega 2,\tfrac \Omega 2]).$$
As
$\lam_{k+1}-\lam_k\geq T$, $k\in\Z$, we conclude that
$\Lambda$ is a set of sampling for $OPW([0,T]{\times}[-\frac \Omega
2,\frac \Omega 2])$ since for $H\in OPW([0,T]{\times}[-\frac \Omega 2,\frac \Omega 2])$,
\begin{eqnarray*}
\|H\|_{HS}^2=\|h_H\|_{L^2}^2&=&\int_0^T\|h_H(t,\cdot)\|^2dt
\ \asymp \ \int_0^T\sum_{n\in\Z}|h_H(t,t+\lam_n)|^2dt\\
&=&\int_0^T\sum_{n\in\Z}|(H\sum_{k\in\Z}\delta_{\lam_k})(t+\lam_n)|^2dt \ \leq\ \|H\sum_{k\in\Z}\delta_{\lam_k}\|_{L^2}^2.
\end{eqnarray*}
\end{proof}

In the remainder of this section, we shall discuss the separation
condition $(\lam_{k+1}-\lam_k)\geq T,\,k\in\Z,$ on the sampling sequence
$\{\lam_k\}_{k\in\Z}$ for\\  $OPW([0,T]{\times}[-\frac \Omega 2,\frac
\Omega 2])$. This condition rules out aliasing in the sense that for
each $x$, the sum $\sum_{n\in\Z}h_H(x-\lam_n,x)= \sum_{n\in\Z}
\kappa_H(x,\lam_n)$
 has only one nonzero summand. Consequently, the sufficient condition on
$\Lambda=\{\lam_k\}_{k\in\Z}$ for being a set of sampling for   $OPW([0,T]{\times}[-\frac \Omega 2,\frac \Omega 2])$  given in Theorem~\ref{thmformerprop}  implies
$$\Omega\leq D^-(\Lambda)\leq D^+(\Lambda)\leq \frac 1 T$$
so that $T\Omega \leq 1$. Theorem~\ref{thmGP} shows that this is not
an additional restriction on  operator Paley-Wiener spaces to allow
for irregular operator sampling.

\begin{cor}\label{CorGP1}
If $\Lambda=\{\lam_k\}_{k\in\Z},\lam_{k+1}>\lam_k$, satisfies
$|\lam_k-kT|\leq L<\frac T 4$ for some $L\geq 0$
and $\lam_{k+1}-\lam_k\geq T$, $k\in\Z$, then $\Lambda$ is a
set of sampling for $OPW([0,T]{\times}[-\frac \Omega 2,\frac \Omega
2])$, $T\Omega \leq 1$.
\end{cor}
\begin{proof}
Since $|\lam_k-kT|\leq L<\frac T 4$, $\{e^{-2\pi i \lam_k\xi}\}$ is
a Riesz basis for $L^2[-\frac 1 {2T},\frac 1 {2T}]$ by Theorem
\ref{Kad}. If $T\Omega <1$, the result follows directly from
Theorem~\ref{thmformerprop}. The case $T\Omega  =1$ follows
analogously.
\end{proof}
Note that the hypothesis on $\Lambda$ in Corollary \ref{CorGP1} is
satisfied if and only if $\lambda_k=kT+\epsilon_k$ with $$-\frac 1
4<-L\leq \ldots \leq \epsilon_{-2}\leq \epsilon_{-1}\leq \epsilon_{0}\leq
\epsilon_{1}\leq \epsilon_{2}\leq \ldots \leq L< \frac 1 4.$$

\begin{thm}\label{setOPW}
Let $\Lambda=\{\lam_k\}_{k\in\Z}$, $\lam_{k+1}> \lam_k$, be a set of
sampling for\\ $OPW([0,T]{\times}[-\frac \Omega 2,\frac \Omega 2])$
and let $\{e^{-2\pi i \lam_k\xi}\}_{k\in\Z}$ be a Riesz basis for
$L^2[-\frac \Omega 2,\frac \Omega 2]$,
   then  $(\lam_{k+1}-\lam_{k}) \geq T$, $k\in\Z$.
\end{thm}
\begin{proof}
It is easy to see that if $\{e^{-2\pi i \lam_k\xi}\}_{k\in\Z}$ is a
Riesz basis for $L^2[-\frac \Omega 2,\frac \Omega 2]$ and
$\sum_{k\in\Z} c_k\delta_{\lam_k}$  is an identifier of
$OPW([0,T]{\times}[-\frac \Omega 2,\frac \Omega 2])$, then $c_k\neq
0$ for all $k\in \Z$.

Assume $\lam_{l+1}-\lam_{l} <T$ for some $l$. For such $l$, set
\begin{eqnarray*}
&\kappa_H(x,\lam_k)&=0\quad \textrm{if}~ k\neq l, l+1\\
\textrm{and}~&\kappa_H(x,\lam_l)&=\left\{\begin{array}{ll}c_{l+1}&\textrm{if}~\lam_{l+1}\leq
x\leq \lam_l+T\\0&\textrm{otherwise}\end{array}\right.\\
\textrm{and}~&\kappa_H(x,\lam_{l+1})&=\left\{\begin{array}{ll}-c_l&\textrm{if}~\lam_{l+1}\leq
x\leq \lam_l+T\\0&\textrm{otherwise}.\end{array}\right.
\end{eqnarray*}
The freedom of choice of values for
$\kappa_{H}(x,\lambda_k)$ is justified by Theorem
\ref{Sei} (2). Then $H\sum_{k\in\Z} c_k\delta_{\lam_k}
(x)=\sum_{k\in\Z} c_k\kappa_H(x,\lam_k) =0$ for all $x\in \Real$,
but as $c_l, c_{l+1}\neq 0$, we have $\kappa_H\neq 0$
and therefore $H\neq 0$.
\end{proof}

\begin{exa}\rm
Theorem~\ref{setOPW}
implies that
$\{\lam_n\}_{n\in\Z}=\{2n\}_{n\in\Z}\cup\{2n+\alpha\}_{n\in\Z}$
with $0<\alpha<1$ is  not a set of sampling
for $ OPW({[0,1]{\times}[-\frac{1}{2},\frac{1}{2}]})$. 
\end{exa}

The condition $(\lam_{k+1}-\lam_k)\geq T,\,k\in\Z$  is not
necessary for operator sampling of
$OPW([0,T]{\times}[-\frac \Omega 2,\frac \Omega 2])$ if $\{e^{-2\pi i \lam_k\xi}\}_{k\in\Z}$ is not a Riesz
basis but a frame for $L^2[-\frac \Omega 2,\frac \Omega 2]$.
\begin{exa}\rm
The tempered distribution $\sum_{k\in\Z} (-1)^k\delta_{\frac k 2}$
identifies\\ $OPW([0,1]\times [-\frac 1 2, \frac 1 2])$ \cite{PfWa06}.
\end{exa}

To illustrate the rigidity of operator sampling in comparison to
function sampling, we add the following simple example.
\begin{exa}\rm
  Let $\Lambda_r=\{\lam_k\}_{k\in\Z}$ be given by $\lam_k=k$ for $k\neq 0$ and $\lam_0=r\in \Real$.
  The set $\Lambda_r$ is a set of sampling for $PW([-\frac 1 2,\frac 1 2])$ if and only if
  $r\notin \Z\setminus\{0\}$. To see this, note that as $\{e^{2\pi i k\xi}\}_{k\in \Z}$ is a
  Riesz basis for $L^2[-\frac 1 2,\frac 1 2]$, so is
  $\{e^{2\pi i k\xi}\}_{k\neq 0} \cup \{ e^{2\pi i r \xi}\}$ if  $r\notin \Z\setminus\{0\}$.
  By Theorem \ref{setOPW}, $\Lambda_r=\{\lam_k\}_{k\in\Z}$ is a set of sampling for $OPW([0,1]{\times}[-\frac 1 2,\frac 1 2])$ if and only if $r=0$.

\end{exa}

\subsection{An operator version of Kramer's Lemma}

Kramer's lemma plays a crucial role in the proofs of a number of important
sampling theorems. For example, it allows for sampling series
expansions for functions which are integral transforms of type other
than Fourier one. For example, Bessel-Hankel, Legendre, Jacobi,
Laguerre, Gegenbauer, Chebyschev, prolate spheroidal, and Hermite
transforms can be considered, where each transform is defined as an
integral transform whose kernel is its special function \cite{Hig,
Meh, Zay}. In particular, if a function on $\R^2$ has a circular
symmetry, then a multi-dimensional Fourier transform can be reduced
to a one-dimensional Bessel-Hankel transform \cite{Pap68}.

\begin{thm}\label{Kra} (Kramer's Lemma) Let $I\subseteq\mathbb{R}$ be a bounded interval and $k(\cdot,t)\in L^2(I)$ for each fixed $t$ in $D\subseteq\mathbb{R}$.
If there is a sampling sequence $\{t_n\}_{n\in\mathbb{Z}}$ in $D$
such that $\{k(\xi,t_n)\}_{n\in\mathbb{Z}}$ forms a frame for
$L^2(I)$, then
for any $f(t)=\langle F,\overline{k(\cdot,t)}\rangle_{L^2(I)}$, $F\in L^2(I)$, 
we have
$$f(t)=\sum_{n\in\mathbb{Z}}f(t_n)S_n(t)$$
where the reconstruction functions $S_n$ are given by
$$S_n(t)=\int_I \tilde{k}_n(\xi)k(\xi,t)d\xi,$$
with $\{\tilde{k}_n(\xi)\}_{n\in\mathbb{Z}}$ a dual frame of
$\{\overline{k(\xi,t_n)}\}_{n\in\mathbb{Z}}$.
\end{thm}
The original Kramer's lemma assumed that $\{k(\xi,t_n)\}_{n\in\Z}$
is an orthonormal basis for $L^2(I)$ \cite{Hig, Kra}. However, we
can easily see that the result extends to the case where
$\{k(\xi,t_n)\}_{n\in\Z}$ forms a frame. We remark that the
classical sampling theorem addressing $f$ in
$PW([-\tfrac{1}{2},\tfrac{1}{2}])$ is given by Kramer's lemma if we
take $k(\xi,t)=e^{2\pi i\xi t}, F=\hat{f}$ and use the fact that
$\{e^{2\pi i n\xi}\}_{n\in\Z}$ is an orthonormal basis for
$L^2[-\tfrac{1}{2},\tfrac{1}{2}]$.

Let $I\subseteq \R$  be bounded. For
\begin{equation} \label{mathcalK}
\mathcal K\,:\,L^2(I) \longrightarrow  L^2(\Real),\quad F \longmapsto\langle F(\cdot),k(x,\cdot)\rangle_{L^2(I)},
\end{equation} bounded, set
\begin{eqnarray*}
  \mathcal{H}^k(S)
 &=&\big\{H\in HS(L^2(\mathbb{R})) \\ &&\quad :\,h_H(t,x)=\langle \zeta_{H}(t,\cdot),k(x,\cdot)\rangle_{L^2},\ \zeta_{H}\in L^2(\R{\times}I),\
\supp\zeta_H\subseteq S\big\}.
\end{eqnarray*} Clearly, for $k(x,\nu)=e^{-2\pi i\nu x}$ we have
$\mathcal{H}^k(S)=OPW(S)$.

\begin{thm}\label{nonuniform}
Let $H\in\mathcal{H}^k({[0,d]\times I})$, $d>0$, $I$ in $\mathbb{R}$
bounded. If there is a set  $\{y_n\}_{n\in\Z}$ in $\mathbb{R}$,
 such that $\{k(t+y_n,\nu)\}_{n\in\Z}$
forms a frame for $L^2(I)$ for every $t\in [0,d]$, and $d'=\inf(y_{n+1}-y_n)\geq d$, then
exists $c>0$ with
\begin{equation}
\|H(\sum_{k\in\Z}\delta_{y_k})\|_{L^2}\geq c \|H\|_{HS},\quad H \in
\mathcal H^k({[0,d]\times I})\label{eqnoneside}.
\end{equation}
If the map $\mathcal K$ in (\ref{mathcalK})
is bounded below, then
\begin{equation}
  \|H(\sum_{k\in\Z}\delta_{y_k})\|_{L^2}\asymp\|H\|_{HS},\quad H \in
\mathcal H^k({[0,d]\times I}), \label{eqnbothsides}
\end{equation}
and operator  reconstruction is possible as
$$h_H(t,x)=r(t)\sum_{n\in\Z}(H\sum_{k\in\Z}\delta_{y_k})(t+y_n)\phi_n(t,x),$$
where $\phi_n(t,x)$ is given by
$$\phi_n(t,x)=\int_I k_n^*(t,\nu)\overline{k(x,\nu)}d\nu$$
with $\{k_n^*(t,\nu)\}_{n\in\Z}$ being a dual frame of
$\{k(t+y_n,\nu)\}_{n\in\Z}$ for each $t$ and $r\equiv 1$ on $[0,d]$
with $\supp r\subseteq [d-d',d']$.
\end{thm}

\begin{proof}
As $\inf_{k}(y_{k+1}-y_k)=d'\geq d$ we have
\begin{equation*}
r(t)(H\sum_{k\in\Z}\delta_{y_k})(t+y_n)=r(t)\sum_{k\in\Z}
h_H(t+y_n-y_k,t+y_n)=h_H(t,t+y_n), \quad t \in\R,
\end{equation*}
where $r\equiv 1$ on $[0,d]$ and $\supp r\subseteq [d-d',d']$.
Consequently,
\begin{eqnarray*}
 \|\zeta_H\|^2_{L^2} &=& \int \|\zeta_H(t,\,\cdot\,)\|^2_{L^2}\,dt \asymp\int \| \{ \langle\zeta_H(t,\, \cdot\,),k(t+y_n,\,\cdot\,)\rangle_{L^2} \} \|^2_{l^2}\,dt \notag \\
  &=& \int \| \{ h_H(t,t+y_n) \} \|^2_{l^2}\, dt
   = \int \sum _{n\in\Z} | h_H(t,t+y_n) |^2 dt  \notag \\
   &=& \int \sum _{n\in\Z} | r(t-y_n)(H\sum_{k\in\Z}\delta_{y_k})(t) |^2 dt
    =  \sum_{n\in\Z}\int_{y_n}^{y_{n+1}}|H\sum_{k\in\Z}\delta_{y_k}(t) |^2 dt \notag \\
      &=      & \int |H\sum_{k\in\Z}\delta_{y_k}(t) |^2 dt  =\|H\sum_{k\in\Z}\delta_{y_k}\|^2, \quad H\in \mathcal H^k([0,d]{\times} I).
\end{eqnarray*}
We used the fact that $r(t-y_n)(H\sum_{k\in\Z}\delta_{y_k})(t)=(H\sum_{k\in\Z}\delta_{y_k})(t)$ for $t\in [y_n,y_n+d]$ and $r(t-y_n)(H\sum_{k\in\Z}\delta_{y_k})(t)=0$ for $t\in [y_n+d,y_{n+1})$, $n\in\Z$.

Equations (\ref{eqnoneside}) and (\ref{eqnbothsides}) follow from
the fact that $\mathcal K$ is bounded and the hypothesis that
$\mathcal K$ is bounded below, respectively.

Moreover, we have, for $\nu\in I$,
\begin{eqnarray*}
\zeta_H(t,\nu)=\sum_{n\in\Z}
\langle\zeta_H(t,\nu),k(t+y_n,\nu)\rangle_{L^2}
k_n^*(t,\nu)=\sum_{n\in\Z} h_H(t,t+y_n)k_n^*(t,\nu),
\end{eqnarray*}
where $\{k_n^*(t,\nu)\}_{n\in\Z}$ is a dual frame of
$\{k(t+y_n,\nu)\}_{n\in\Z}$ for each $t$.
 Multiplying $\overline{k(x,\nu)}$ and integrating with respect
to $\nu$ on both sides, we have for fixed $t\in [0,d]$
$$h_H(t,x)=\sum_{n\in\Z} h_H(t,t+y_n)\phi_n(t,x),\quad x\in\Real,$$
where $\phi_n(t,x)=\int_I k_n^*(t,\nu)\overline{k(x,\nu)}d\nu$.
\end{proof}

Note that a set of sampling $\{y_n\}_{n\in\Z}$ in Theorem
\ref{nonuniform} is uniformly discrete and a separation constant
$d'$ is greater than or equal to $d$. On the other hand, the
condition that $\{k(t+y_n,\nu)\}_{n\in\Z}$ forms a frame for $L^2(I)$ indicates
that the sampling points should not be too sparse.
\begin{rem}\rm
Unlike Kramer's Lemma for functions, we do not have any explicit
example for operator sampling other than the $OPW(S)$ case. Sampling theorems based on  various kinds of orthogonal
polynomials generally do not satisfy all hypotheses in Theorem
\ref{nonuniform}. For instance, for the Bessel-Hankel transform, we
have $k(x,\nu)=\sqrt{x\nu}J_n(x\nu)$ where $J_n$ is the Bessel
function of the first kind of order $n$. Taking $\lam_k$ as the
$k$-th positive root of $J_n(x)$, one can, in fact,  obtain a sampling
expansion induced from Bessel-Hankel transform \cite{Hig}. However,
the kernel $k(x,\nu)=\sqrt{x\nu}J_n(x\nu)$ does not allow for
$\{k(t+y_n,\nu)\}_{n\in\Z}$ being a frame for each $t$.
\end{rem}

In Theorem \ref{uniform} we used the fact that $\{m(\nu)e^{-2\pi
in\nu}\}_{n\in\Z}$ forms an orthonormal basis for
$L^2[-\tfrac{1}{2},\tfrac{1}{2}]$ for any function $m$ satisfying
$|m(\nu)|=1$.  In general, we have the following.

\begin{prop}\label{lemma}
For $D, I\subseteq \Real$ and $k(x,\nu)$ defined on $\Real\times I$,
let $k(x,\nu)$ satisfy
$k(t+x,\nu)=m(t,\nu)k(x,\nu)$ 
for some $m$ such that
$0<\|m(t,\cdot)\|_0\leq\|m(t,\cdot)\|_\infty<\infty$ for all $t\in
D$. If $\{k(y_n,\nu)\}_{n\in\Z}$ is a frame for $L^2(I)$, then
$\{k(t+y_n,\nu)\}_{n\in\Z}$ is also a frame for $L^2(I)$ for all
$t\in D$.
\end{prop}
For example, the Fourier kernel $k(x,\nu)=e^{-2\pi ix\nu}$ and the
Hilbert transform kernel $k(x,\nu)=-i\sgn(\nu)\,e^{-2\pi ix\nu}$
satisfy the hypotheses of Proposition \ref{lemma}.

\section{Multi-channel sampling for Hilbert-Schmidt operators}

In classical multi-channel sampling, a signal is reconstructed using discrete
values from the outputs of $N$ different time--invariant operators applied to a single input signal. Generally, each of the $N$ outputs is sampled at
$\tfrac{1}{N}$-th of the Nyquist-Landau rate of the input signal. For example, when
the signal's bandwidth is $\Omega$ Hz, then we should collect at
least $\Omega$ samples per second. But if we design $N$ channel
filters appropriately, then it suffices to obtain $\Omega/N$ samples per second from
each channel and combine the samples in order to reconstruct the
signal.  This allows us to reduce sampling rate requirements on sampling hardware at the cost of employing multiple samplers.

A number of important theorems, for example, on periodic nonuniform sampling, on
derivative sampling and on samples of Hilbert transforms can be
explained in the framework of  multi-channel sampling \cite{Bro,Hig}. In \cite{HKK},  multi-channel sampling has been
developed for abstract Hilbert space, allowing  each channel output
to be sampled at different, irregular points.

In Sections \ref{section:uniform} and \ref{section:irregular}, it
has been shown that a slowly time-varying/underspread operator is
identifiable by a single channel output while in \cite{KoPf,PfWa} it
is shown that an overspread operator is not identifiable in this sense. However,
in this section we shall show that overspread operators may be
recoverable from multiple channel outputs. In addition, we seek to
reduce the rate at which delta impulses are produced for channel
identification.

Throughout this section, we shall consider  Hilbert-Schmidt operators whose Kohn-Nirenberg symbols
are bandlimited to rectangular domains.

\begin{thm}\label{multichannel}
For  $M,N\in\mathbb{N}$, $OPW([0,N]{\times}[-\tfrac{M}{2},\tfrac{M}{2}])$ permits multi-channel operator sampling as
\begin{equation*}
\|H\|_{HS}^2=\frac{1}{M^2N}\sum_{j=0}^{MN-1}\|H(\sum_{n\in\mathbb{Z}}{e^{2\pi
ijn/MN}}\delta_{\tfrac{n}{M}})\|^2,\quad H\in OPW([0,N]{\times}[-\tfrac{M}{2},\tfrac{M}{2}]).
\end{equation*}
\end{thm}
\begin{proof}
 Suppose that $\eta_H\in
L^2([0,N]{\times}[-\frac{M}{2},\frac{M}{2}])$ for some
$M,N\in\mathbb{N}$. Consider an orthonormal basis
$\{\frac{1}{\sqrt{M}}\,e^{-2\pi i(t+\tfrac{n}{M})\nu}\}_{n\in\Z}$
for $L^2[-\frac{M}{2},\frac{M}{2}]$, $t\in[0,N]$. Then
\begin{eqnarray*}
\eta_H(t,\nu)&=&\sum_{n\in\Z}\langle \eta_H(t,\nu),
\tfrac{1}{\sqrt{M}}\,e^{-2\pi i(t+\tfrac{n}{M})\nu}\rangle
\tfrac{1}{\sqrt{M}}\,e^{-2\pi i(t+\tfrac{n}{M})\nu}\\
&=&\sum_{n\in\Z} h_H(t,t+\tfrac{n}{M})\tfrac{1}{{M}}\,e^{-2\pi
i(t+\tfrac{n}{M})\nu},\quad \nu\in [-\tfrac{M}{2},\tfrac{M}{2}],
\end{eqnarray*}
and
\begin{eqnarray*}
h_H(t,x)=\sum_{n\in\Z} h_H(t,t+\tfrac{n}{M})~\sinc\,
M(x-t-\tfrac{n}{M}),~ x\in\Real,~a.e.~ t\in [0,N].
\end{eqnarray*}
Since $\{\sqrt{M}\,\sinc M(\cdot-t-\tfrac{n}{M})\}_{n\in\Z}$ is an
orthonormal basis for $PW([-\frac{M}{2},\frac{M}{2}])$,  Parseval's
identity gives
$$\|h_H(t,\cdot)\|_2^2=\tfrac{1}{M}\sum_{n\in\Z} |h_H(t,t+\tfrac{n}{M})|^2,\quad  t\in
[0,N].$$ Hence, we have
\begin{eqnarray*}
\|h_H\|_{L^2(\mathbb{R}^2)}^2&=&\int_0^N\|h_H(t,\cdot)\|^2_{L^2(\mathbb{R})}dt\
= \ \tfrac{1}{M}\int_0^N\sum_{n\in\Z}|h_H(t,t+\tfrac{n}{M})|^2 dt\\
&=&\tfrac{1}{M}\sum_{r=0}^{MN-1}\int_{\tfrac{r}{M}}^{\tfrac{r+1}{M}}\sum_{n\in\Z}|h_H(t,t+\tfrac{n}{M})|^2
dt\\
&=&\tfrac{1}{M}\int_0^{\frac{1}{M}}\sum_{r=0}^{MN-1}\sum_{n\in\Z}|h_H(t+\tfrac{r}{M},t+\tfrac{r}{M}+\tfrac{n}{M})|^2
dt\\
&=&\tfrac{1}{M}\int_0^{\frac{1}{M}}\sum_{r=0}^{MN-1}\sum_{n\in\Z}|h_H(t+\tfrac{r}{M},t+\tfrac{n}{M})|^2
dt.
\end{eqnarray*}
Since
$(H\sum_n\frac{1}{\sqrt{MN}}\delta_{\frac{n}{M}})(t+\frac{k}{M})=
\sum_n\frac{1}{\sqrt{MN}}h_H(t+\frac{k}{M}-\frac{n}{M},t+\frac{k}{M})$,
 for fixed $t\in[0,\frac{1}{M}]$, we have
$$r(t)\sum_{n\in\Z}\tfrac{1}{\sqrt{MN}}h_H(t+\tfrac{k}{M}-\tfrac{n}{M},t+\tfrac{k}{M})=
\tfrac{1}{\sqrt{MN}}\sum_{r=0}^{MN-1}h_H(t+\tfrac{r}{M},t+\tfrac{k}{M}),$$
where $r(t)=\chi_{[0,\frac{1}{M}]}(t)$. Similarly, consider
$(H\sum_n\frac{1}{\sqrt{MN}}\omega_j^n\delta_{\frac{n}{M}})(x)$
where $\omega_j=e^{2\pi i j/MN}$, $j=0,1,\cdots,MN-1$. Then
$$(H\sum_n\tfrac{1}{\sqrt{MN}}\omega_j^n\delta_{\frac{n}{M}})(t+\tfrac{k}{M})=\sum_n
\tfrac{\omega_j^n}{\sqrt{MN}}h_H(t+\tfrac{k}{M}-\tfrac{n}{M},t+\tfrac{k}{M}).$$
For fixed $t\in[0,\frac{1}{M}]$,
$$r(t)\sum_{n\in\Z}\tfrac{\omega_j^n}{\sqrt{MN}}h_H(t+\tfrac{k}{M}-\tfrac{n}{M},t+\tfrac{k}{M})=
\sum_{r=0}^{MN-1}\tfrac{\omega_j^{k-r}}{\sqrt{MN}}h_H(t+\tfrac{r}{M},t+\tfrac{k}{M}).$$
Consider the system of linear equations
\begin{eqnarray*}
r(t)\left[\begin{array}{c}(H\sum_n\frac{1}{\sqrt{MN}}\delta_{\frac{n}{M}})(t+\frac{k}{M})\\
(H\sum_n\frac{\omega_1^n}{\sqrt{MN}}\delta_{\frac{n}{M}})(t+\frac{k}{M})\\
\vdots\\
(H\sum_n\frac{\omega_{MN-1}^n}{\sqrt{MN}}\delta_{\frac{n}{M}})(t+\frac{k}{M})
\end{array}\right]=
A_k\left[ \begin{array}{c}h_H(t,t+\frac{k}{M})\\
h_H(t+\frac{1}{M},t+\frac{k}{M})\\
\vdots\\
h_H(t+\frac{MN-1}{M},t+\frac{k}{M})
 \end{array}\right]
\end{eqnarray*}
where all matrices $A_k$ are unitary $MN\times MN$ DFT matrices
with entries given by
$(A_k)_{j,l}=\frac{1}{\sqrt{MN}}\,e^{2\pi i(j-1)(k-(l-1))/MN}$.
Since all $A_k$'s are unitary, we have
\begin{eqnarray*}\label{systemofeqn}
\left\|\left[\begin{array}{c}\|\{(H\sum_n\frac{1}{\sqrt{MN}}\delta_{\frac{n}{M}})(t{+}\frac{k}{M})\}_k\|_{l^2}\\
\|\{(H\sum_n\frac{\omega_1^n}{\sqrt{MN}}\delta_{\frac{n}{M}})(t{+}\frac{k}{M})\}_k\|_{l^2}\\
\vdots\\
\|\{(H\sum_n\frac{\omega_{MN-1}^n}{\sqrt{MN}}\delta_{\frac{n}{M}})(t{+}\frac{k}{M})\}_k\|_{l^2}\\
\end{array}\right]\right\|^2
{=}
\left\| \left[\begin{array}{c}
\|\{h_H(t,t{+}\frac{k}{M})\}_k\|_{l^2}\\
\|\{h_H(t{+}\frac{1}{M},t{+}\frac{k}{M})\}_k\|_{l^2}\\
\vdots\\
\|\{h_H(t{+}\frac{MN-1}{M},t{+}\frac{k}{M})\}_k\|_{l^2}
 \end{array} \right]
\right\|^2,
\end{eqnarray*}
where $\| \cdot \|$ denotes the $L^2[0,\frac{1}{M}]^{MN}$ norm.
Therefore
\begin{eqnarray*}
\|h_H\|^2&=&\tfrac{1}{M}\int_0^{\tfrac{1}{M}}\sum_{k\in\Z}\sum_{j=0}^{MN-1}|(H\sum_{n\in\Z}
\tfrac{\omega_j^n}{\sqrt{MN}}\, \delta_{\frac{n}{M}})
(t+\tfrac k M)|^2dt\\
&=&\tfrac{1}{M}\int_\mathbb{R}\sum_{j=0}^{MN-1}|(H\sum_{n\in\Z}\tfrac{\omega_j^n}{\sqrt{MN}}\, \delta_{\frac{n}{M}})
(t)|^2dt\\&=&\tfrac{1}{M^2N}\sum_{j=0}^{MN-1}\|H\sum_{n\in\Z}{\omega_j^n}\delta_{\frac{n}{M}}\|_{L^2(\mathbb{R})}^2.
\end{eqnarray*}
\end{proof}

Clearly, the matrices $A_k$ can be replaced by appropriate sequences of
matrices whose norms are bounded above and away from zero.
\begin{thm}\label{multicor}
Let  $M,N\in\mathbb{N}$
and $\{c_{j,n}\}_{j=1,n\in\Z}^{MN}$ bounded with $A_k$, $k\in\Z$, invertible with $ \|A_k^{-1}\| \leq  C<\infty$ where $(A_k)_{j,l=1,\ldots,MN}=c_{j,k-l+1}$. For
$f_j=\sum_{n\in\Z}c_{j,n}\delta_{\frac{n}{M}}$, $1\leq j\leq MN$, we have
$$\|H\|_{HS}^2\asymp \sum_{j=1}^{MN}\|Hf_j\|_{L^2}^2,\quad  H\in OPW([0,N]{\times}[-\tfrac{M}{2},\tfrac{M}{2}]).$$
\end{thm}
The assumption in Theorem \ref{multicor} is satisfied if for all
$j$, $\{c_{j,n}\}_{n\in\Z}$ is $MN-$periodic and $A_{j,0}$ is
invertible with $\|A_{j,0}^{-1}\|$ also uniformly bounded.

Now we consider periodic nonuniform
sampling as first proposed by Yen
\cite{Yen}. We
first recall periodic nonuniform sampling theorem for functions in
$PW([-\frac \Omega 2,\frac \Omega 2])$.
\begin{thm}\label{PNS}
There exists a
Riesz basis $\{S_j(t-\frac{nN}{\Omega})\}_{j=1,n\in\Z}^N$ for
$PW([-\frac \Omega 2,\frac \Omega 2])$ such that
$$f(t)=\sum_{j=1}^N\sum_{n\in\Z}f\big(\tfrac{nN}{\Omega}+\alpha_j\big)S_j\big(t-\tfrac{nN}{\Omega}\big), \quad f\in PW([-\tfrac \Omega 2,\tfrac \Omega 2]),$$
where $0\leq \alpha_j<\frac N \Omega$, $1\leq j\leq N$, and
$\alpha_i\neq\alpha_j$ for $i\neq j$.
\end{thm}
We show that $OPW( [0,N]{\times}[-\tfrac{M}{2},\tfrac{M}{2}])$ is
identifiable by $MN$ identifiers which are given by delta-trains
whose supports are periodically nonuniformly distributed.
\begin{thm}\label{theorem:nonuniformly}
For
$M,N\in\mathbb{N}$, and $0\leq\alpha_1< \alpha_2<\ldots<\alpha_{MN}<N$, we have
\begin{equation*}
\|H\|_{HS}^2\asymp\sum_{j=1}^{MN}\|H(\sum_{n\in\mathbb{Z}}\delta_{nN+\alpha_j})\|^2, \quad H\in OPW([0,N]{\times}[-\tfrac{M}{2},\tfrac{M}{2}]).
\end{equation*}
\end{thm}
\begin{proof}
If we apply Theorem \ref{PNS} to $h_H(t,\cdot)\in PW([-\frac M
2,\frac M 2])$ with $MN$ channels, then we obtain
$$h_H(t,x)=\sum_{j=1}^{MN}\sum_{n\in\Z}h_H(t,t+nN+\alpha_j)\,\varphi_j(x-t-nN),~ x\in\Real,~a.e.~ t\in [0,N],$$
and
$$\|h_H(t,\cdot)\|^2\asymp\sum_{j=1}^{MN}\sum_{n\in\Z}|h_H(t,t+nN+\alpha_j)|^2,~a.e.~ t\in [0,N],$$
where $\{\varphi_j(x-t-nN)\}_{j=1,n\in\Z}^{MN}$ is a Riesz basis for
$PW([-\frac M 2,\frac M 2])$ for each $t\in [0,N]$. Therefore, we
have
\begin{eqnarray*}
\|H\|_{HS}^2=\|h_H\|_{L^2}^2&\asymp&\int_0^N\sum_{j=1}^{MN}\sum_{n\in\Z}|h_H(t,t+nN+\alpha_j)|^2dt\\
&=&\sum_{j=1}^{MN}\int_0^N\sum_{n\in\Z}|H(\sum_{k\in\mathbb{Z}}\delta_{kN+\alpha_j})(t+nN+\alpha_j)|^2dt\\
&=&\sum_{j=1}^{MN}\|H(\sum_{k\in\mathbb{Z}}\delta_{kN+\alpha_j})\|_{L^2}^2,
\end{eqnarray*}
that is, $\{\sum_{k\in\mathbb{Z}}\delta_{kN+\alpha_j}\}_{j=1}^{MN}$
identifies $OPW([0,N]{\times}[-\tfrac{M}{2},\tfrac{M}{2}])$.
\end{proof}

\begin{rem}\rm
  Note that as in the multichannel sampling theory for functions, Theorem~\ref{theorem:nonuniformly} can be applied to reduce the rate of impulse transmission. For example, as $OPW(
[0,1]{\times}[-\tfrac{1}{2},\tfrac{1}{2}])\subseteq OPW(
[0,N]{\times}[-\tfrac{1}{2},\tfrac{1}{2}])$, we can identify any
operator in $OPW( [0,1]{\times}[-\tfrac{1}{2},\tfrac{1}{2}])$ by its
action on the tempered distributions
$\sum_{n\in\Z}\delta_{nN+j\alpha }$, $0\leq \alpha\leq
1$,\,$j=0,\ldots, N-1$, each of which has impulse rate $1/N$.
\end{rem}

Another important multi-channel sampling concept can be applied to operator
sampling, namely, derivative sampling.
We apply a multi-channel sampling formula
consisting of samples of $f$ and its $MN-1$ derivatives \cite{Lin, LiAb}.
\begin{thm}
For $M,N\in\mathbb{N}$, we have
\begin{equation*}
\|H\|_{HS}^2\asymp\sum_{j=0}^{MN-1}\Big\|\sum_{r=0}^j {j\choose
r}(-1)^r\big(H\sum_k\delta_{kN}^{(r)}\big)^{(j-r)}\Big\|^2,~ H\in
OPW([0,N]{\times}[-\tfrac{M}{2},\tfrac{M}{2}]).
\end{equation*}
Here, $f^{(r)}$ denotes the $r$-th derivative of $f$ in the distributional sense.
\end{thm}

\begin{proof}
If we apply multi-channel derivative sampling theorem  to $h_H(t,\cdot)\in PW([-\frac M
2,\frac M 2])$, then we obtain
$$h_H(t,x)=\sum_{j=0}^{MN-1}\sum_{n\in\Z}\frac{\partial^j}{\partial x^j} h_H(t,x)|_{x=t+nN}\,\varphi_j(x-t-nN),~ x\in\Real,~a.e.~ t\in [0,N],$$
and
$$\|h_H(t,\cdot)\|^2\asymp\sum_{j=0}^{MN-1}\sum_{n\in\Z}\big|\frac{\partial^j}{\partial x^j} h_H(t,x)|_{x=t+nN}\big|^2,~a.e.~ t\in [0,N],$$
where $\{\varphi_j(x-t-nN)\}_{j=0,n\in\Z}^{MN-1}$ is a Riesz basis
for $PW([-\frac M 2,\frac M 2])$ for each fixed $t\in [0,N]$. Let
$H_j$ be the operator defined through the $j$-th derivative of
$h_H(t,\cdot)$, that is,
$$H_jf(x)=\int \frac{\partial^j}{\partial x^j} h_H(t,x) f(x-t)dt\quad a.e.$$
Since $\supp \frac{\partial^j}{\partial x^j} h_H(\cdot,x)\subseteq
\supp h_H(\cdot,x)$, we have
$$\chi_{[0,N]}(t)(H_j\sum_{k\in\Z}\delta_{kN})(t+nN)=\frac{\partial^j}{\partial
x^j} h_H(t,x)|_{x=t+nN},$$
so that
\begin{eqnarray*}
\|h_H\|^2&=&\int_0^N\|h_H(t,\cdot)\|^2dt \asymp \int_0^N\sum_{j=0}^{MN-1}\sum_{n\in\Z}|(H_j\sum_{k\in\Z}\delta_{kN})(t+nN)|^2dt\\
&=&\sum_{j=0}^{MN-1}\|H_j\sum_{k\in\Z}\delta_{kN}\|^2.
\end{eqnarray*}
Observing that $H_jf(x)=\sum_{r=0}^j{j\choose
r}(-1)^r\left(Hf^{(r)}\right)^{(j-r)}$ by Leibniz's rule completes the proof. 
\end{proof}

We conclude this section by presenting two explicit multi-channel
reconstruction formulas for operators.

\begin{exa} \rm
Consider the Riesz
bases $\{e^{-2\pi
i(t+2n)\nu}\}_{n\in\Z}\cup\{-2\pi i\nu\, e^{-2\pi
i(t+2n)\nu}\}_{n\in\Z}$, $t\in\R$, of $L^2[-1/2,1/2]$,  as  well as their Riesz basis duals given by
$\{2(1-2|\nu|)\,e^{-2\pi
i(t+2n)\nu}\}_{n\in\Z}\cup\{\frac{2i}{\pi}~\mathrm{sgn}(\nu)\,e^{-2\pi
i(t+2n)\nu}\}_{n\in\Z}$ \cite{Hig}. We obtain
\begin{eqnarray*}
\eta_H(t,\nu)&=&\sum_n\langle \eta_H(t,\cdot),e^{-2\pi
i(t+2n)\cdot}\rangle 2(1-2|\nu|)\,e^{-2\pi i(t+2n)\nu}\\
&+&\sum_n\langle \eta_H(t,\cdot),-2\pi i\cdot e^{-2\pi
i(t+2n)\cdot}\rangle\, \tfrac{2i}{\pi}\,\mathrm{sgn}(\nu)\,e^{-2\pi
i(t+2n)\nu},~\nu\in [-\tfrac 1 2, \tfrac 1 2].
\end{eqnarray*}
Taking an inverse Fourier transform with respect to the variable
$\nu$, we have
$$h_H(t,x)=\sum_n h_H(t,t+2n)S_n(t,x)+\dfrac{\partial}{\partial x}h_H(t,x)|_{x=t+2n}T_n(t,x),~x\in\Real,\, t\in [0,2].$$
where $$S_n(t,x)=\sinc^2\tfrac{1}{2}(x-t-2n)$$ and
$$T_n(t,x)=\tfrac{2}{\pi}\,\sinc\tfrac{1}{2}(x-t-2n)\sin\tfrac{\pi}{2}(x-t-2n).$$
\end{exa}
We give a second reconstruction formula for the  operator class
$OPW({[0,2]{\times}[-\frac{1}{2},\frac{1}{2}]})$.
\begin{exa}\rm
Let $H\in OPW({[0,2]{\times}[-\frac{1}{2},\frac{1}{2}]})$. Then for
$0<\alpha<1$,
$$h_H(t,x)=\sum_{n\in\Z} (H\sum_{k\in\mathbb{Z}}\delta_{2k}) (t+2n)S_1(x-t-2n)+
(H\sum_{k\in\mathbb{Z}}\delta_{2k+\alpha})(t+2n+\alpha)S_2(x-t-2n-\alpha),$$
where the reconstruction functions are
\begin{equation*}
S_1(x)=\frac{2}{e^{\pi i\alpha}-1}\mathcal{F}^{-1}\left(e^{\pi
i\alpha}\chi_{[-\frac{1}{2},0)}(\nu)-\chi_{[0,\frac{1}{2}]}(\nu)\right)(x)
\end{equation*}
and
\begin{equation*}
 S_2(x)=\frac{2}{e^{\pi
i\alpha}-1}\mathcal{F}^{-1}\left(-\chi_{[-\frac{1}{2},0)}(\nu)+e^{\pi
i\alpha}\chi_{[0,\frac{1}{2}]}(\nu)\right)(x).
\end{equation*}
\end{exa}

\bibliographystyle{elsart-num}

\end{document}